\documentclass{article}
\usepackage{amssymb,amsmath,amsthm,graphicx,kotex} 
\usepackage{scalerel}
\usepackage{xcolor}

\def\qed{\hfill {\hbox{${\vcenter{\vbox{               
   \hrule height 0.4pt\hbox{\vrule width 0.4pt height 6pt
   \kern5pt\vrule width 0.4pt}\hrule height 0.4pt}}}$}}}

\newtheorem{theorem}{Theorem}
\newtheorem{lemma}[theorem]{Lemma}
\newtheorem{proposition}[theorem]{Proposition}

\theoremstyle{definition}
\newtheorem{example}{Example}
\newtheorem{definition}{Definition}

\date{}

\title{\Large \textbf{Mosaics for immersed surface-links}} 

\author{
Seonmi Choi\footnote{Email: smchoi@knu.ac.kr.}
 \and
Jieon Kim\footnote{Email: jieonkim7@gmail.com}
}

\begin{document}
\maketitle

\begin{abstract}
The concept of a knot mosaic was introduced by Lomonaco and Kauffman as a means to construct a quantum knot system. The mosaic number of a given knot $K$ is defined as the minimum integer $n$ that allows the representation of $K$ on an $n \times n$ mosaic board. Building upon this, the first author and Nelson extended the knot mosaic system to encompass surface-links through the utilization of marked graph diagrams and established both lower and upper bounds for the mosaic number of the surface-links presented in Yoshikawa's table. 
In this paper, we establish a mosaic system for immersed surface-links by using singular marked graph diagrams. We also provide the definition and discussion on the mosaic number for immersed surface-links.
\end{abstract}

\parbox{6in} {\textsc{Keywords:} mosaic knot, surface-link, immersed surface-link, marked graph diagram, singular vertex

\smallskip

\textsc{2020 MSC:} 57K12, 57K45, 81P99}

\section{\large\textbf{Introduction}}\label{I}

A knot mosaic is a mathematical object that represents a knot diagram composed of tiles, where each tile represents a segment of the knot. Knot mosaic diagrams provide a powerful tool for understanding and analyzing the structure of quantum knots. They allow for a visual representation of the knot and its components, enabling researchers to study the intricate connections and relationships within the knot. 
By utilizing knot mosaic diagrams, various properties and characteristics specific to quantum knots can be defined and investigated. 
Thus, knot mosaic diagrams serve as an important framework for the definition and study of quantum knots.

In 2008, Lomonaco and Kauffman \cite{LK} introduced a knot mosaic as a construction for quantum knot systems and conjectured that tame knot theory is equivalent to knot mosaic theory. 
In 2014, Kuriya and Shehab \cite{KurShe} proved this conjecture by using grid diagrams, also known as rectangular diagrams. 
The mosaic number of a knot $K$ is defined as the smallest number $n$ for which $K$ can be represented as an $n$-mosaic. 
As a result, mathematicians have extensively studied the calculation of the mosaic number for knots, focusing on different types of knots and exploring both upper and lower bounds, among other research directions.
See \cite{KurShe, HLLO, LHLO, LLPP, OHLL} for more information. 

A {\it surface-link} is the image of a closed surface smoothly (piecewise linearly and locally flatly) embedded in $\mathbb{R}^4$ (or $S^4$). If the underlying surface is connected, it is referred to as a {\it surface-knot}. 
Two surface-links $F$ and $F'$ are said to be {\it equivalent} if there exists an orientation-preserving homeomorphism $h : \mathbb{R}^4 \rightarrow \mathbb{R}^4$ such that $h(F)=F'$.
In order to facilitate the study of surface-links, several descriptive tools have been developed due to the inherent complexity of directly dealing with them in 4-space. 
These tools include broken surface diagrams, marked graph diagrams, banded link diagrams, motion pictures, and others. 
Using these alternative representations allows for effective investigation and analysis of surface-links in 4-space. 
For more details, refer to \cite{CKS00, Kamadabook, Kamadabook2, KawShiSuz, Yoshikawa}.
A {\it marked graph} is a $4$-valent graph in $\mathbb{R}^3$ with a {\it marker}, which is a line segment at a $4$-valent vertex. 
A {\it marked graph diagram} is a generic projection on $\mathbb{R}^2$ with classical crossings and marked vertices as a link diagram. 
It is well known that surface-links can be represented by marked graph diagrams (cf. \cite{CKS00, Kamadabook, Kamadabook2, KeaKur, Swenton, Yoshikawa}).

In 2022, the first author and Nelson \cite{ChoNel} constructed a mosaic system for surface-links using marked graph diagrams as a tool. 
They established lower and upper bounds of the mosaic number for the surface-links in Yoshikawa’s table and enhanced the kei counting invariant for unoriented surface-links as well as classical knots and links using mosaic diagrams, as an application.

An {\it immersed surface-link} is a closed surface smoothly immersed in $\mathbb{R}^4$ (or $S^4$) with the property of self-transversality.
An immersed surface is a fundamental role in the field of low-dimensional topology, particularly in the study of $4$-manifolds. 
Despite their significance, it is challenging to provide explicit descriptions of immersed surface-links and their equivalence, except for a few specific examples.
Therefore, there is a pressing need for developing effective methods to represent immersed surface-links.
On the other hand, extensive investigations have been conducted on representations for embedded surface-links, and recently, there has been considerable interest in extending these representations for immersed surface-links.

In 2018, Kamada and Kawamura \cite{KamKaw} extended the concept of a normal form from surface-links to immersed surface-links and demonstrated that any immersed surface-link can be represented in a normal form.
Kamada, Kawauchi, Kim, and Lee \cite{KKKL18} introduced a diagrammatic method for representing immersed surface-links in $\mathbb{R}^4$ via marked graph diagrams, using the concept of $H$-admissibility, and then Kim \cite{Kim20} enumerated immersed surface-links up to ch-index $9$.
In 2019, the second author and Kawauchi \cite{KawKim19} gave infinitely many immersed surface-links with certain types of double point singularities. 
In 2020, Hughes, Kim, and Miller \cite{HKM20} introduced a banded unlink diagrams in Kirby diagrams for surfaces embedded in any $4$-manifold and
in 2021, they extended their results to studying surfaces immersed in any $4$-manifolds via a singular banded unlink diagram \cite{HKM21} and this is another diagrammatic system compared to the results in \cite{KKKL18}.
In 2022, Jab\l onowski \cite{Jab22} constructed marked graph diagrams with singular vertices, called {\it singular marked graph diagrams}, for representing surfaces immersed in $\mathbb{R}^{4}$ from their results. 

In this paper, we extend mosaics for embedded surface-links to immersed surface-links by using singular marked graph diagrams. 
The paper is organized as follows. 
In Section \ref{SMGD}, we review the relationship between immersed surface-links and singular marked graph diagrams. 
Section \ref{MGM} introduced a mosaic system for embedded surface-links. 
Mosaics for immersed surface-links and their equivalence can be defined in Section \ref{SMGM}. 
In Section \ref{Mnbr}, we calculate the mosaic number of the immersed surface-links, illustrated in \cite{Kim20}, and suggest their lower and upper bounds.


\section{\large\textbf{Immersed surface-links and singular marked graph diagrams}}\label{SMGD}

We introduce the definitions and notions of singular marked graph diagrams for representing immersed surface-links in $4$-space. See \cite{HKM21, Jab22} in details. 

An {\it immersed surface-link} is the image of a closed surface smoothly immersed in $\mathbb{R}^4$ (or $S^4$) such that the multiple points are transverse double points. 
An {\it immersed surface-knot} is an immersed surface-link whose underlying surface is connected. 
Clearly, if it has no transverse double points, then it is a surface-link.
Two immersed surface-links $F$ and $F'$ are said to be {\it equivalent} in $\mathbb{R}^4$ if there exists an orientation-preserving homeomorphism $h : \mathbb{R}^4 \rightarrow \mathbb{R}^4$ such that $h(F)=F'$.

We use an effective tool for describing immersed surface-links, known as a \textit{singular marked 
graph diagram}. 
A \textit{singular marked graph} is a spatial graph embedded in 
$\mathbb{R}^3$ possibly with $4$-valent vertices decorated by a line segment, such as {\Large \scalerel*{\includegraphics{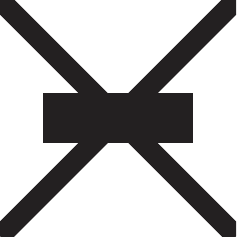}}{B}}, 
or an encircled decoration, such as 
{\Large \scalerel*{\includegraphics{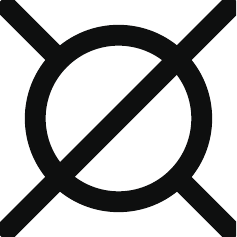}}{B}}.
We call such a line segment a \textit{marker} and a vertex with a marker a 
\textit{marked vertex}. A vertex with an encircled decoration is called a \textit{singular vertex}.
In the same way as a link diagram, one can define a \textit{singular marked graph 
diagram} which is a diagram in $\mathbb{R}^2$ with classical crossings, marked vertices, and singular vertices.

\begin{figure}[h!]
  \centering
  \includegraphics[width = 12cm]{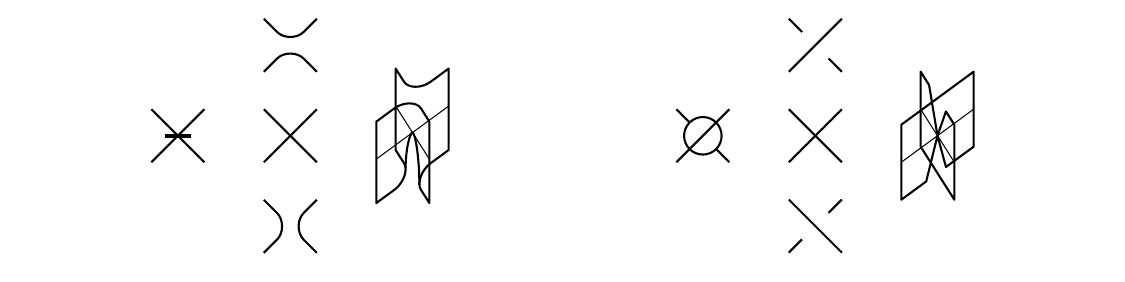}
  \caption{Marked vertex and singular vertex}\label{MarkedSingular}
\end{figure}

Let $D$ be a singular marked graph diagram.
In the case of a marked vertex 
{\Large \scalerel*{\includegraphics{MarkedVertex1.pdf}}{B} }, 
the local diagram obtained by splicing in a direction consistent with its marker (referred to as the $+$ direction), appears as 
{\Large \scalerel*{\includegraphics{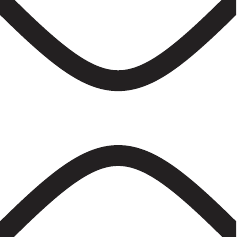}}{B} }. 
Conversely, by applying this process in the opposite direction (referred to as the $-$ direction), the resulting local diagram appears as 
{\Large \scalerel*{\includegraphics{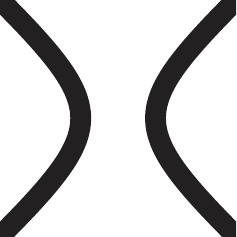}}{B} }. 
In the case of a singular vertex 
{\Large \scalerel*{\includegraphics{SingularVertex1.pdf}}{B} },
the local diagram obtained by choosing the arc in the circle symbol as the over arc (referred to as the $+$ direction), appears as 
{\Large \scalerel*{\includegraphics{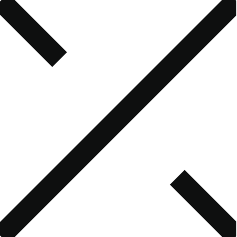}}{B} }. 
Conversely, the resulting local diagram by connecting the two arcs disconnected by the circle symbol (referred to as the $-$ direction), appears as 
{\Large \scalerel*{\includegraphics{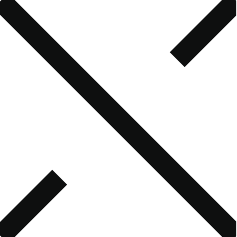}}{B} }. 
Consequently, two distinct classical link diagrams, denoted by $L_{+}(D)$ and $L_{-}(D)$, can be obtained from $D$ by applying all marked and singular vertices in the $+$ and $-$ directions, respectively. 
Indeed, their local diagrams and neighborhoods are described in Figure \ref{MarkedSingular}.
We refer to $L_{+}(D)$ and $L_{-}(D)$ as the \textit{positive} and \textit{negative resolutions} of $D$, respectively.

\begin{figure}[h!]
  \centering
  \includegraphics[width = 8cm]{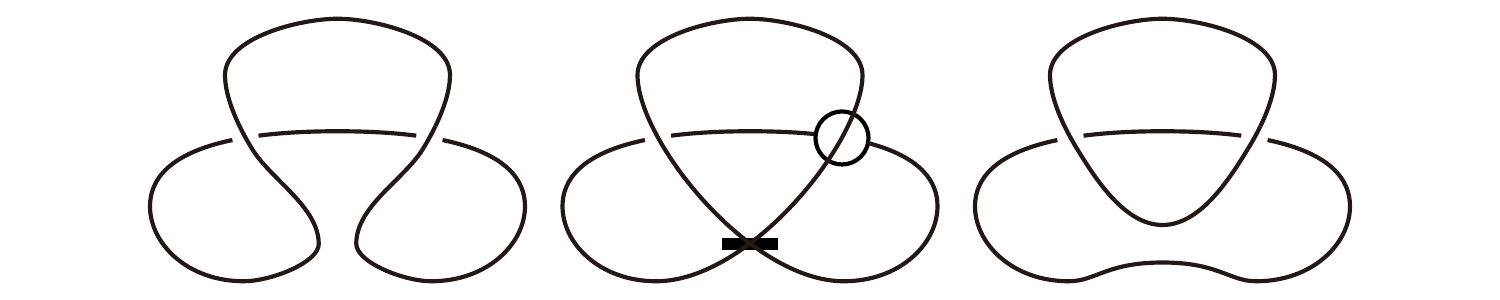}
  \caption{An example of singular marked graph diagrams}\label{ExaSMGD}
\end{figure}

A singular marked graph diagram $D$ is said to be \textit{admissible} if both resolutions $L_{-}(D)$ and $L_{+}(D)$ are diagrams of trivial links. 
A singular marked graph is said to be \textit{admissible} if it has an admissible singular marked graph diagram.
For example, the singular marked graph diagram $D$ as depicted in Figure \ref{ExaSMGD} is admissible.

\begin{figure}[h!]
  \centering
  \includegraphics[width = 8cm]{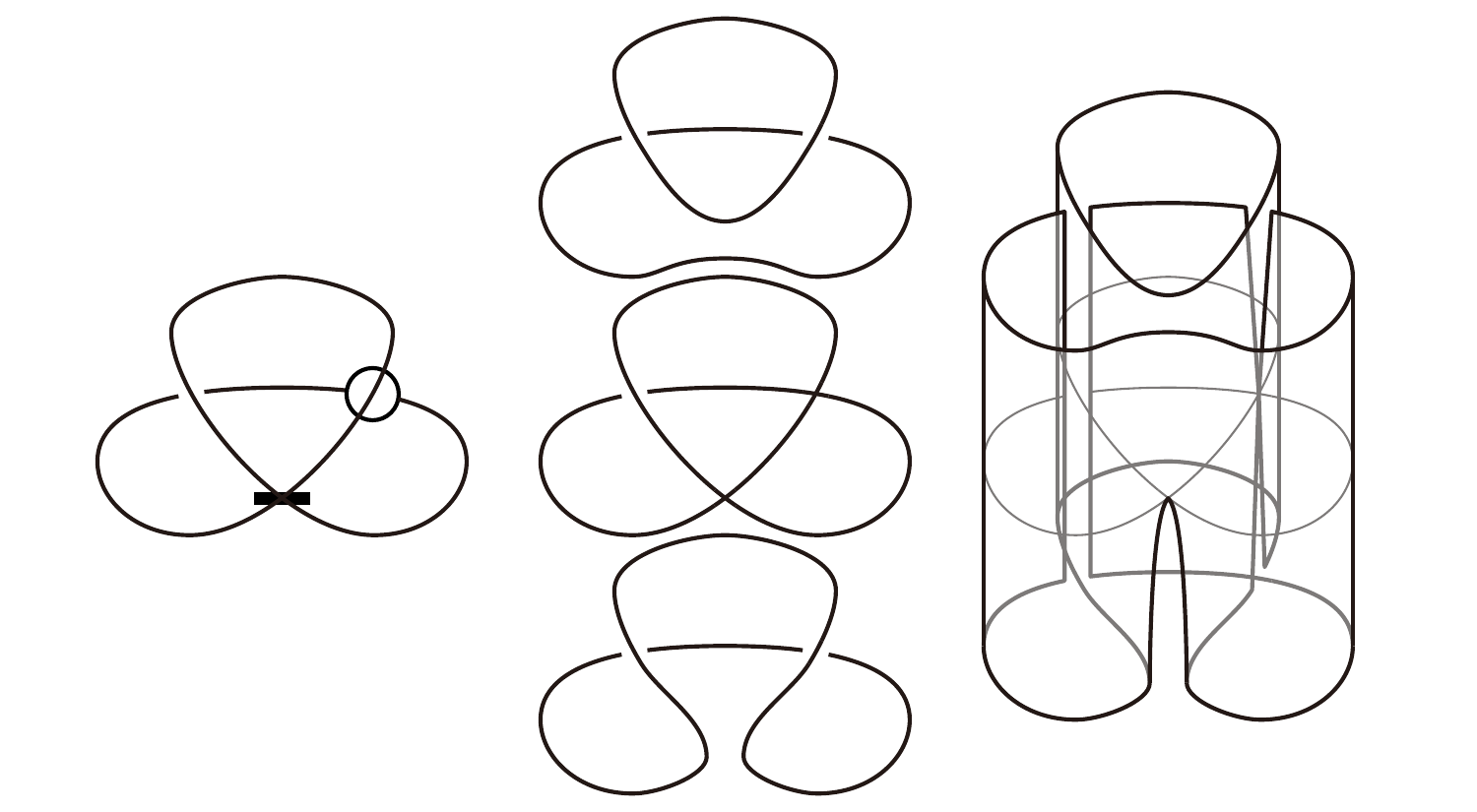}
  \caption{Cobordism obtained from a singular marked graph diagram}\label{Cobordism}
\end{figure}

For a given admissible singular marked graph diagram $D$,
one can construct a uniquely immersed surface-link $F(D)$ obtained from $D$, up to equivalence. Conversely, every immersed surface-link $F$ can be expressed by an admissible singular marked graph diagram $D$, that is, $F(D)$ is equivalent to $F$.
See \cite{HKM21, Jab22} for more details. 
For example, let $D$ denote the singular marked graph diagram with one classical crossing, one marked vertex, and one singular vertex as described in Figure \ref{Cobordism}.
Then one can obtain its cobordism constructed by connecting the positive and negative resolutions of $D$. 

\begin{figure}[h!]
  \centering
  \includegraphics[width = 14cm]{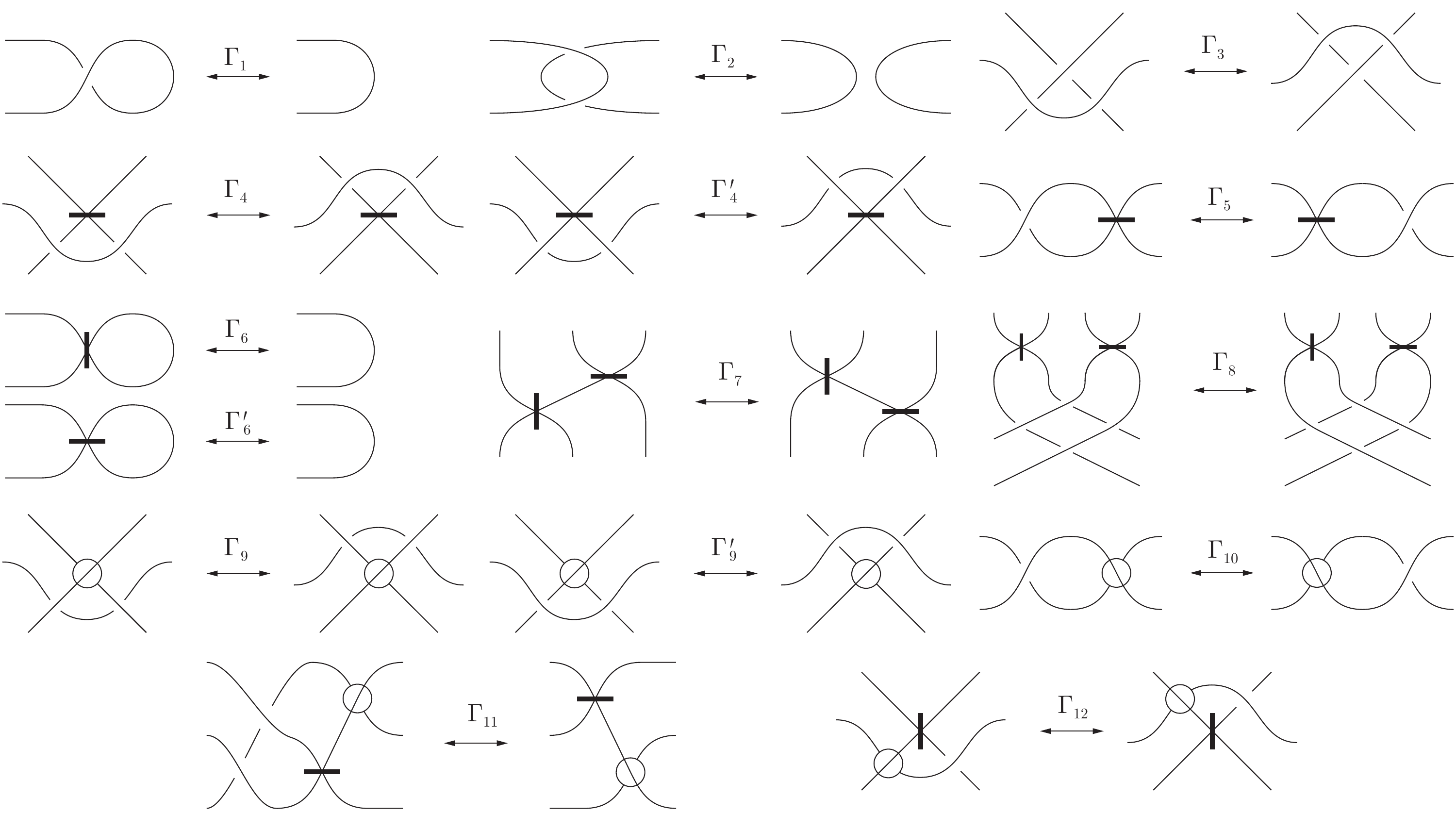}
  \caption{Singular Yoshikawa moves for immersed surface-links}\label{GenYoshikawaMoves}
\end{figure}

In \cite{Jab22}, the equivalence moves for immersed surface-links are the generalized Yoshikawa moves of singular marked graph diagrams and its minimal generating set is depicted in Figure \ref{GenYoshikawaMoves}.
In this paper, we call the generalized Yoshikawa moves the {\it singular Yoshikawa moves.}

\begin{proposition}[\cite{HKM21, Jab22}]
Two singular marked graph diagrams $D$ and $D'$ present equivalent immersed surface-links if and only if $D$ can be obtained from $D'$ by a finite sequence of ambient isotopies in $\mathbb{R}^{2}$ and singular Yoshikawa moves.
\end{proposition}

\begin{figure}[h!]
  \centering
  \includegraphics[width = 12cm]{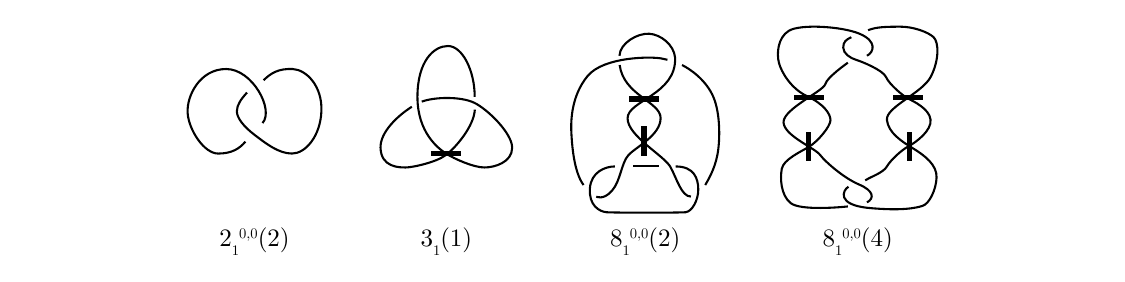}
  \caption{Examples via marked graph diagrams}\label{KimTableSing}
\end{figure}

In \cite{Kim20}, the concept of H-admissibility is used in a marked graph diagram for immersed surface-links to categorize cases up to ch-index 9
and the examples in the table of immersed surface-links that has at least one singular vertex can be illustrated in Figure \ref{KimTableSing}.
A marked graph diagram is {\it H-admissible} if both positive and negative resolutions are H-trivial link diagrams where an H-trivial diagram is a diagram of a split union of trivial knots and Hopf links. 
During the process of constructing an immersed surface-link from an H-admissible marked graph diagram, a singular point arises when capping the Hopf link. 
See \cite{KKKL18, Kim20} in details.

\begin{figure}[h!]
  \centering
  \includegraphics[width = 10cm]{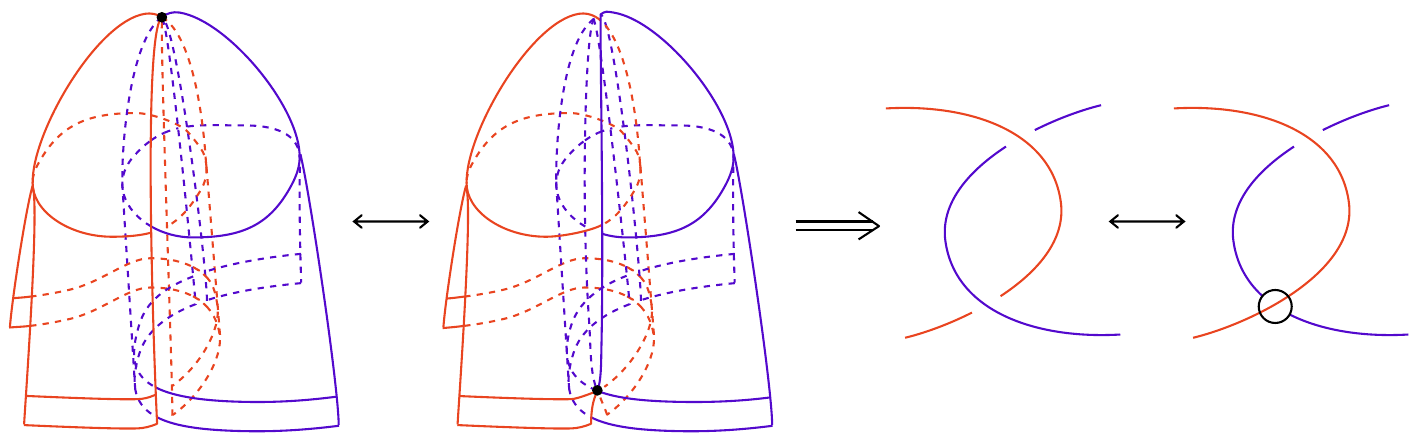}
  \caption{Upper singularity}\label{UpSing}
\end{figure}

\begin{figure}[h!]
  \centering
  \includegraphics[width = 10cm]{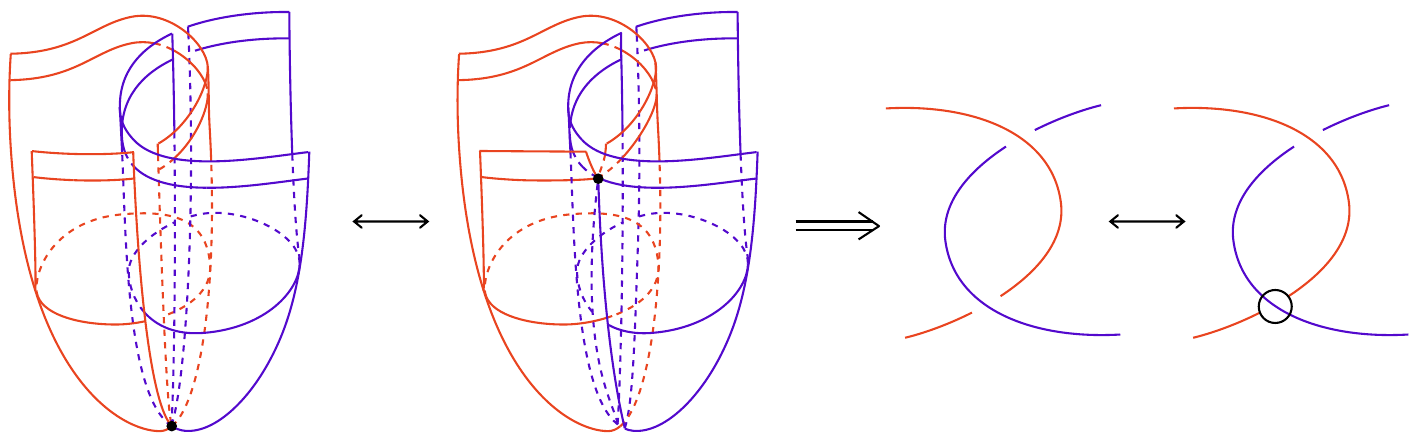}
  \caption{Lower singularity}\label{LowSing}
\end{figure}

\begin{figure}[h!]
  \centering
  \includegraphics[width = 10cm]{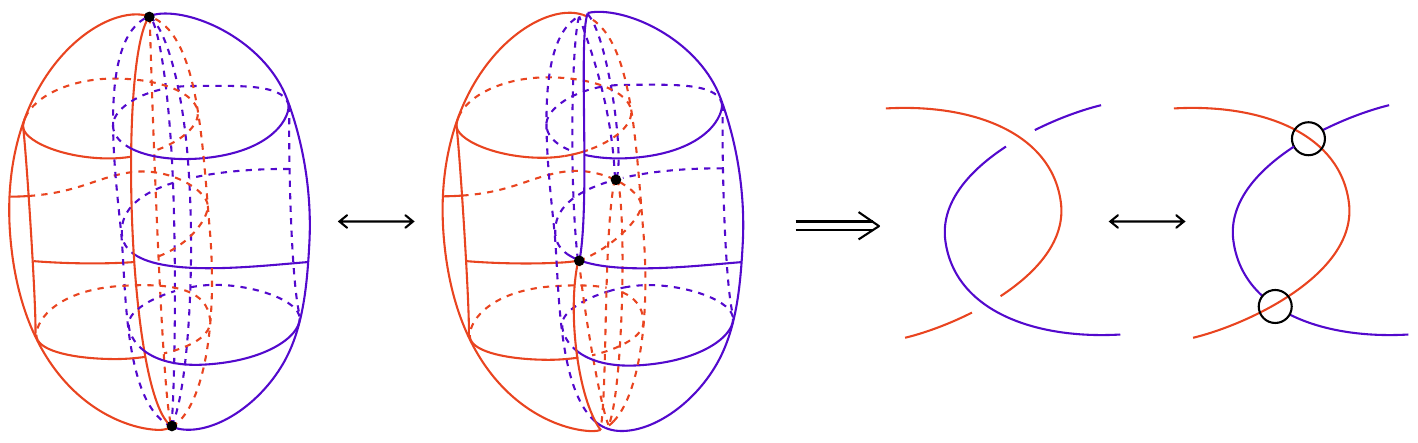}
  \caption{Upper and lower singularities}\label{BothSing}
\end{figure}

\begin{figure}[h!]
  \centering
  \includegraphics[width = 12cm]{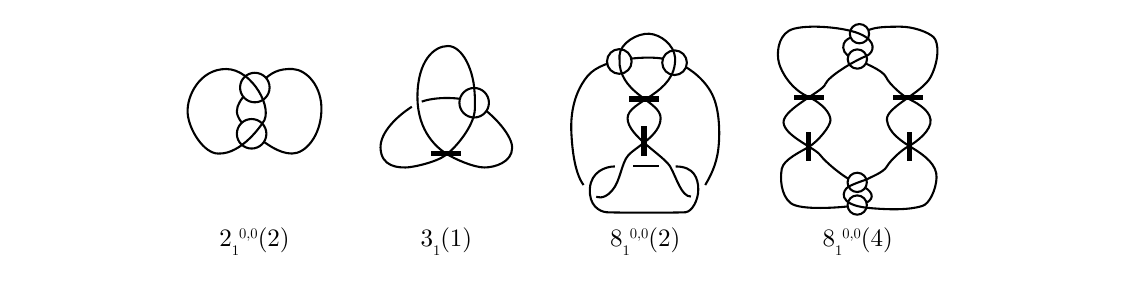}
  \caption{Examples via singular marked graph diagrams}\label{KimTableSing2}
\end{figure}

In order to convert H-admissible marked graph diagrams into admissible singular marked graph diagrams, we consider the singularities of an H-admissible marked graph diagram.
Its singularities can be divided into three cases based on their location and their broken surface diagrams are shown in the first figure of each of Figures \ref{UpSing}, \ref{LowSing}, and \ref{BothSing}. 
Figure \ref{UpSing} is the case of a singularity appearing from the positive resolution, called the {\it upper singularity}. 
The upper singularity can be moved down to the level where the marked graph diagram appears, to converting this H-admissible marked graph diagram to a singular marked graph diagram and its local move used to transform between them is also shown in Figure \ref{UpSing}. 
Conversely, the case of a singularity appearing from the negative resolution, called the {\it lower singularity}, is illustrated in Figure \ref{LowSing}, similar to Figure \ref{UpSing}.
In the last case, both upper and lower singularities are present, as depicted in Figure \ref{BothSing}, and the transformation consists of two processes. 
Therefore, one can deform the marked graph diagrams of Figure \ref{KimTableSing} into the singular marked graph diagrams as depicted in Figure \ref{KimTableSing2}.





\section{\large\textbf{Mosaics for surface-links via marked graph diagram}}\label{MGM}

We review the mosaic system for surface-links using marked graph diagrams and focus on the case of unoriented surface-links. The consideration of orientation can be deferred to the last stage.

\begin{figure}[h!]
  \centering
  \includegraphics[width = 12cm]{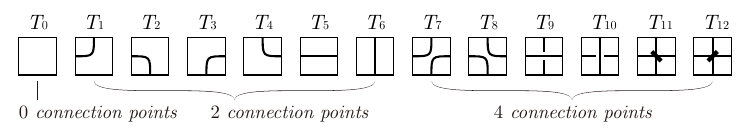}
  \caption{Mosaic tiles}\label{MosaicTiles}
\end{figure}

Let $\mathbb{T}$ be the set of {\it mosaic tiles}, which are rectangles $T_{0}, T_{1}, \cdots, T_{12}$ with arcs and possibly with one crossing or marked vertex as depicted in Figure \ref{MosaicTiles}.
Indeed, there are exactly $6$ tiles, up to rotation.
The endpoints of an arc on a mosaic tile is located the center of an edge, called {\it connection points}, and then each tile has $0$, $2$, or $4$ connection points. 

An {\it $(m, n)$-mosaic} is an $m\times n$ matrix whose entries are mosaic tiles in $\mathbb{T}$.
If it is square, i.e. $m=n$, then it is called an {\it $n$-mosaic}.
The sets consisting of $(m, n)$-mosaics and $n$-mosaics are denoted by $\mathbb{M}^{(m, n)}_{M}$ and $\mathbb{M}^{(n)}_{M}$, respectively.
\begin{figure}[h!]
  \centering
  \includegraphics[width = 12cm]{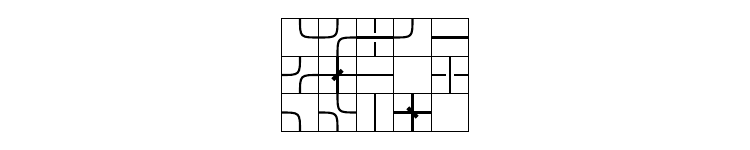}
  \caption{A $(3,5)$-mosaic diagram}\label{Exa1}
\end{figure}
Two tiles in a mosaic are said to be {\it contiguous} if they are adjacent in either the same row or column.
A tile in a mosaic is said to be {\it suitably connected} if all its connection points coincide with the connection points of contiguous tiles. 
Note that for a $(3,5)$-mosaic illustrated in Figure \ref{Exa1}, the tile at entry $(2, 2)$ is suitably connected, whereas the tile at entry $(1,4)$ is not suitably connected. 

\begin{definition}
A {\it marked graph $(m, n)$-mosaic} is an $(m, n)$-mosaic in which all tiles are suitably connected. The set of all knot $(m, n)$-mosaic is the subset of $\mathbb{M}^{(m, n)}_{M}$, denoted by $\mathbb{K}^{(m, n)}_{M}$. 
If $m=n$, then it is called a {\it marked graph $n$-mosaic} and its set is denoted by $\mathbb{K}^{(n)}_{M}$.
\end{definition}

\begin{example}
The spun trefoil $8_{1}$ has a marked graph $8$-mosaic and $7$-mosaic diagrams in Figure \ref{Exa2}.
\begin{figure}[h!]
  \centering
  \includegraphics[width = 12cm]{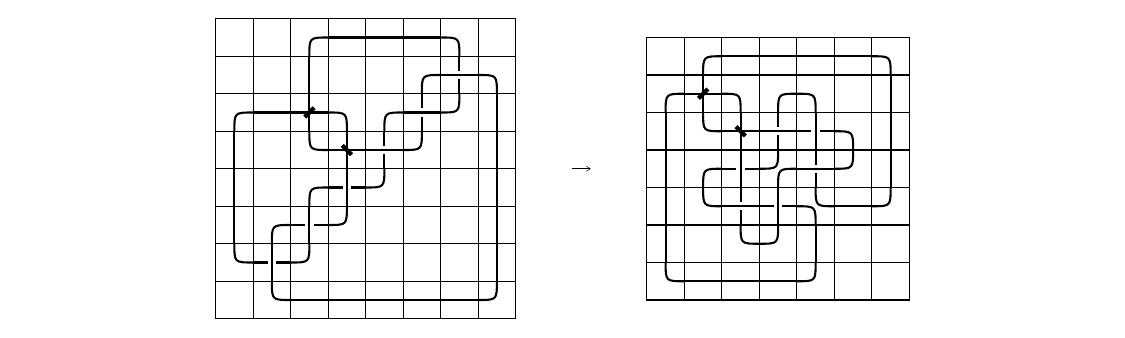}
  \caption{Two mosaic diagrams of $8_{1}$}\label{Exa2}
\end{figure}
\end{example}

For representing surface-links using marked graph diagrams, one can consider the equivalence moves for surface-links, planar moves and Yoshikawa moves, \cite{Yoshikawa}.
The non-determistic tiles are necessary to define the moves and see \cite{LK} in detail.

\begin{figure}[h!]
  \centering
  \includegraphics[width = 12cm]{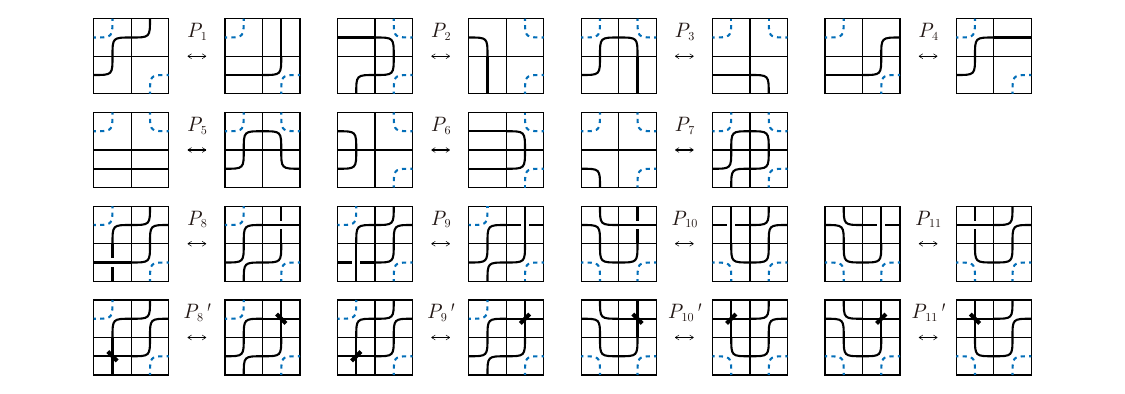}
  \caption{Planar moves for mosaics}\label{PlanarM}
\end{figure}

\begin{figure}[h!]
  \centering
  \includegraphics[width = 12cm]{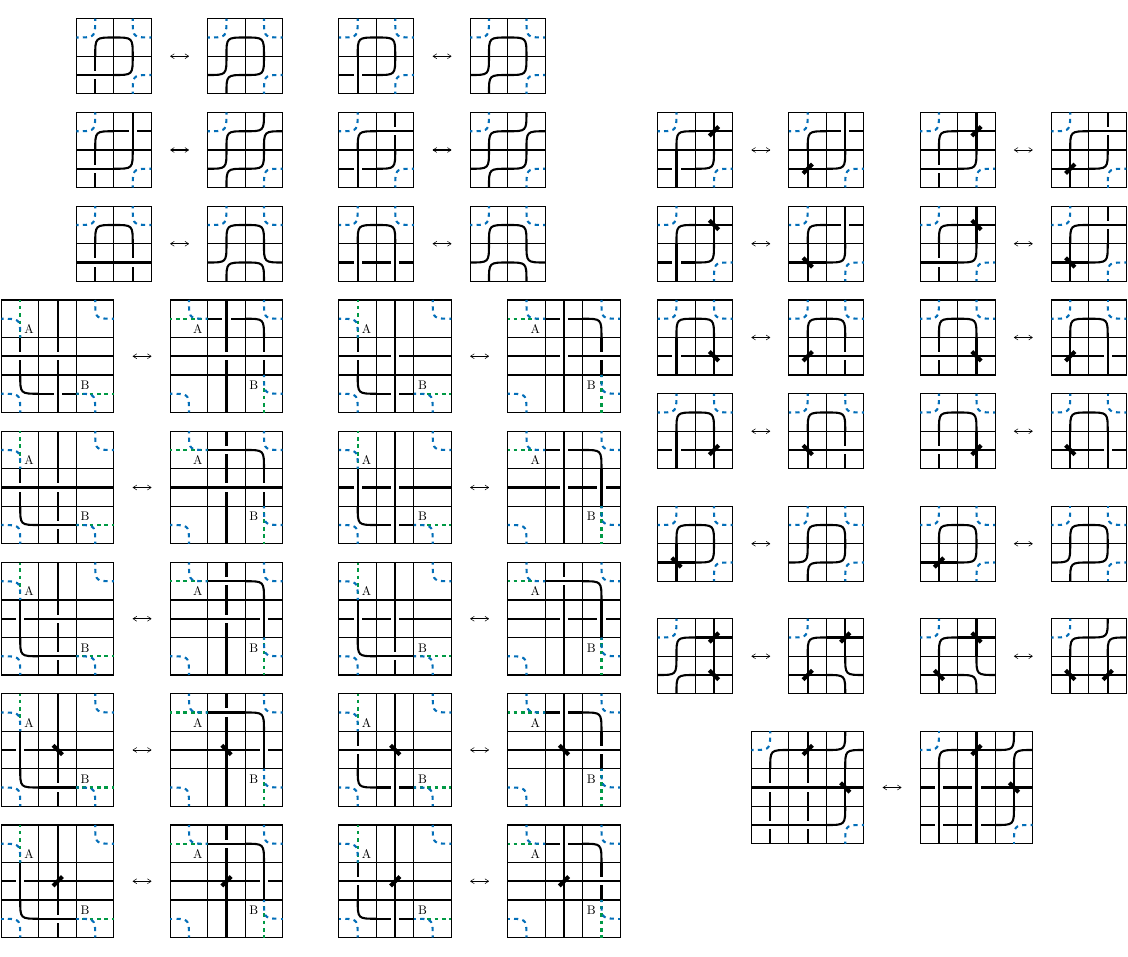}
  \caption{Yoshikawa moves for mosaics}\label{YoshikawaM}
\end{figure}

The equivalence of marked graph mosaics consists of 15 types for planar isotopy, as shown in Figure \ref{PlanarM},
 and Yoshikawa moves illustrated in Figure \ref{YoshikawaM}.

All marked graph mosaic moves are permutations on the set $\mathbb{M}^{(n)}_{M}$ of $n$-mosaics. Indeed, they are also in the group of all permutations of the set $\mathbb{K}^{(n)}_{M}$ of marked graph $n$-mosaics.

\begin{definition}
The \textit{ambient isotopy group $\mathbb{A}^{(n)}_{M}$} is the subgroup of the group of all permutations of the set $\mathbb{K}^{(n)}_{M}$ generated by all planar isotopy moves and all Yoshikawa moves.  
\end{definition}

Two marked graph $n$-mosaics $M$ and $M'$ are said to be \textit{of the same marked graph $n$-type}, denoted by 
  $M \overset{n}{\sim} M',$
  if there exists an element of $\mathbb{A}^{(n)}_{M}$ such that it transforms $M$ into $M'$.
Two marked graph $n$-mosaics $M$ and $M'$ are said to be \textit{of the same marked graph type} 
  if there exists a non-negative integer $k$ such that 
  $$i^{k}M \overset{n+k}{\sim} i^{k}M',$$
  where $i : \mathbb{M}^{(j)} \rightarrow \mathbb{M}^{(j+1)}$ is the mosaic injection by adding a row and a column consisting of only empty tiles.
Therefore, we can obtain the following result.  
  
\begin{proposition}[\cite{ChoNel}]
Let $M$ and $M'$ be two marked graph mosaics of two marked graphs $K$ and $K'$, respectively. 
Then $M$ and $M'$ are of the same marked graph mosaic type if and only if $K$ and $K'$ are equivalent.
\end{proposition}

\begin{example}
The mosaic number of $7_{1}^{0,-2}$ is six whose mosaic diagram is described in Figure \ref{ExaMnbr1}.
\begin{figure}[h!]
  \centering
  \includegraphics[width = 12cm]{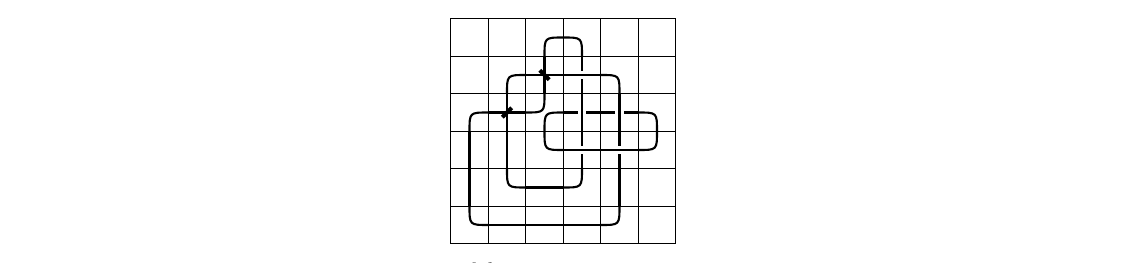}
  \caption{A mosaic diagram of $7_{1}^{0,-2}$}\label{ExaMnbr1}
\end{figure}
\end{example}


\section{\large\textbf{Mosaics for immersed surface-links via singular marked graph diagrams}}\label{SMGM}

We introduce the mosaic system for representing immersed surface-links in $\mathbb{R}^{4}$ using singular marked graph diagrams and construct additional mosaic moves to the mosaic moves as shown in Figures \ref{PlanarM} and \ref{YoshikawaM}, for the equivalence of immersed surface-links.

\begin{figure}[h!]
  \centering
  \includegraphics[width = 12cm]{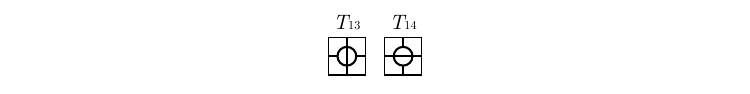}
  \caption{Mosaic tiles with a singular vertex}\label{MosaicS}
\end{figure}

Let $T_{13}$ and $T_{14}$ denote the mosaic tiles with one singular vertex as described in Figure \ref{MosaicS}.
Note that $T_{13}$ and $T_{14}$ are identical by rotation and possess four connection points. 
The set $\mathbb{T}_{S}$ is defined as the collection of all tiles $T_{0}, T_{1}, \cdots, T_{14}$, i.e. $\mathbb{T}_{S}=\mathbb{T}_{M}\cup \{T_{13}, T_{14}\}.$
The introduction of these two additional mosaic tiles to $\mathbb{T}_{M}$ naturally extends the various concepts defined in $\mathbb{T}_{M}$ into $\mathbb{T}_{S}$.

From now on, the entries of an $(m, n)$-mosaic are tiles in $\mathbb{T}_{S}$ and 
the sets consist of $(m, n)$-mosaics and $n$-mosaics are denoted by $\mathbb{M}_{S}^{(m, n)}$ and $\mathbb{M}_{S}^{(n)}$, respectively.

\begin{definition}
A {\it singular marked graph $(m, n)$-mosaic} is an $(m, n)$-mosaic in which all tiles are suitably connected.
If $m=n$, then it is called a {\it singular marked graph $n$-mosaic}.
Their sets are denoted by $\mathbb{K}_{S}^{(m, n)}$ and $\mathbb{K}_{S}^{(n)}$, respectively.
\end{definition}

\begin{example}
The immersed surface-link $3_{1}(1)$ is represented by the singular marked graph $4$-mosaic as illustrated in Figure \ref{Exa3_1(1)}.
  \begin{figure}[h!]
    \centering
    \includegraphics[width = 12cm]{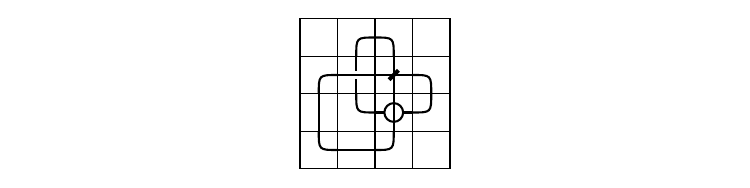}
    \caption{A singular marked graph mosaic}\label{Exa3_1(1)}
  \end{figure}  
\end{example}

\begin{figure}[h!]
  \centering
  \includegraphics[width = 10cm]{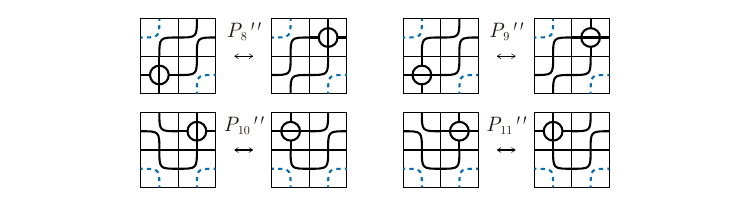}
  \caption{Additional planar mosaic moves for singular vertices}\label{PlanarM_S}
\end{figure}

For the equivalence of immersed surafce-links, there are planar isotopy and singular Yoshikawa moves as depicted in Figrue \ref{GenYoshikawaMoves}.
In order to construct mosaic moves for the equivalence of singular marked graph mosaics corresponding to the equivalence of immersed surafce-links, in the case of planar isotopy, we need four mosaic moves $P_{8}'', P_{9}'', P_{10}''$ and $P_{11}''$ as illustrated in Figure \ref{PlanarM_S}.

\begin{figure}[h!]
  \centering
  \includegraphics[width = 12cm]{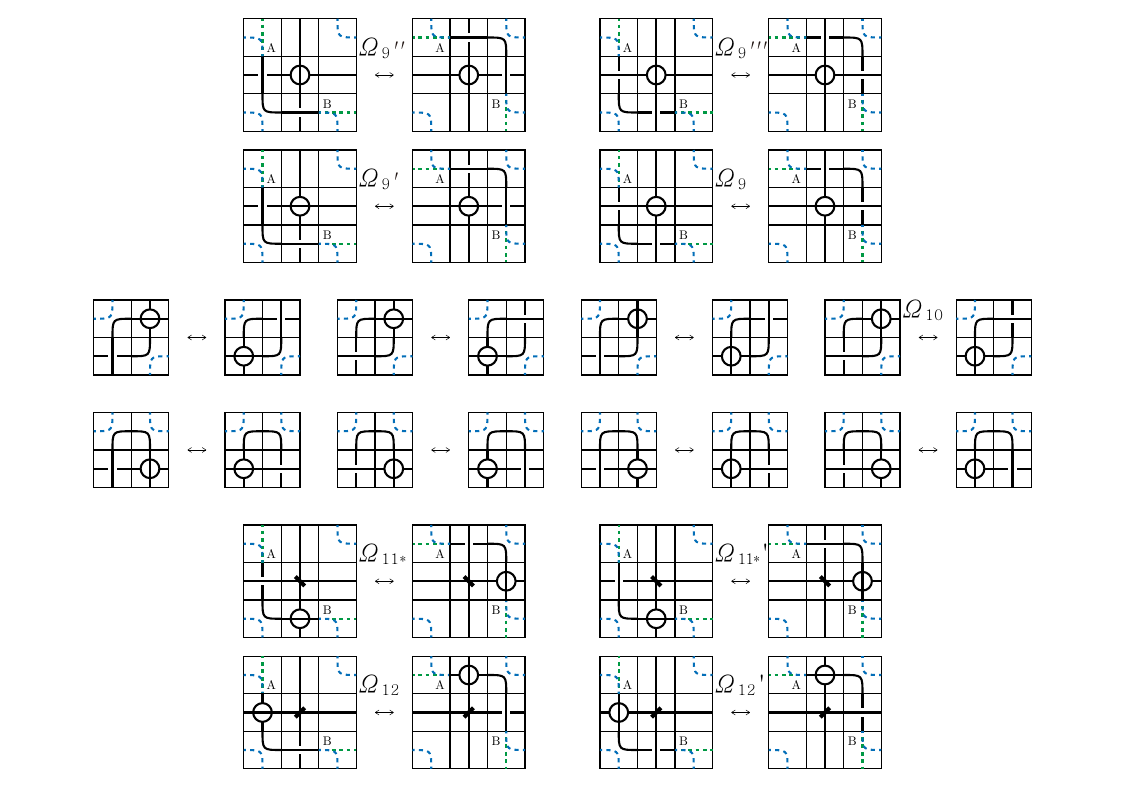}
  \caption{Mosaic moves for Yoshikawa moves with singular vertices}\label{YoshikawaM_S}
\end{figure}

The mosaic moves for singular Yoshikawa moves $\Gamma_{1}, \Gamma_{2}, \cdots, \Gamma_{8}$ are the same with marked graph mosaic moves in Figure \ref{YoshikawaM}.
The mosaic moves for moves $\Gamma_{9}, \Gamma_{9}', \Gamma_{10}, \Gamma_{11}, \Gamma_{12}$ are as described in Figure \ref{YoshikawaM_S}. 

All mosaic moves are permutations on the set $\mathbb{M}^{(n)}_{S}$ of $n$-mosaics. Indeed, they are also in the group of all permutations of the set $\mathbb{K}^{(n)}_{S}$ of singular marked graph $n$-mosaics.

\begin{definition}
The {\it ambient isotopy group $\mathbb{A}^{(n)}_{S}$} is the subgroup of the group of all permutations of the set $\mathbb{K}^{(n)}_{S}$ generated by all planar isotopy moves and all singular Yoshiakawa moves.  
\end{definition}

Two singular marked graph $n$-mosaics $M$ and $M'$ are said to be {\it of the same singular marked graph $n$-type}, denoted by 
  $M \overset{n}{\sim} M',$
  if there exists an element of $\mathbb{A}^{(n)}_{S}$ such that it transforms $M$ into $M'$.
Two singular marked graph $n$-mosaics $M$ and $M'$ are said to be {\it of the same singular marked graph type}
  if there exists a non-negative integer $k$ such that 
  $$i^{k}M \overset{n+k}{\sim} i^{k}M',$$
  where $i : \mathbb{M}^{(j)}_{S} \rightarrow \mathbb{M}^{(j+1)}_{S}$ is the mosaic injection by adding a row and a column consisting of only empty tiles.\\
  
By the construction of mosaic moves for singular marekd graph mosaics, we have

\begin{theorem}\label{ThmMain}
Let $M$ and $M'$ be two singular marked graph mosaics of two immersed surface-links $F$ and $F'$, respectively. 
Then $M$ and $M'$ are of the same singular marked graph type if and only if $F$ and $F'$ are equivalent.
\end{theorem}

Therefore, by Theorem \ref{ThmMain}, two singular marked graph mosaics $M$ and $M'$ are said to be {\it of the same immersed surface-link type}  
if they are of the same singular marked graph type.\\

\begin{figure}[h!]
  \centering
  \includegraphics[width = 12cm]{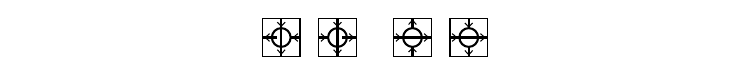}
  \caption{The orientation for mosaics with a singular vertex}\label{MosaicOri}
\end{figure}

To illustrate the oriented case, one can explore all the potential types for assigning orientations to each unoriented mosaic tile. 
The oriented tiles with one singular vertex are shown in Figure \ref{MosaicOri} and the orientation of other tiles was introduced in \cite{ChoNel, LK}. 
Consequently, the establishment of an oriented mosaic theory guarantees the congruence of orientations with the original theory, thereby ensuring seamless alignment between the two approaches.

\section{\large\textbf{Mosaic numbers}}\label{Mnbr}

The mosaic number for an immersed surface-link can be defined in the same manner using singular marked graph diagrams. 
In this section, we only consider the case of immersed surface-links with at least one singular vertex.

\begin{definition}
The {\it mosaic number} of an immersed surface-link $F$, denoted by $m(F)$, is the smallest integer $n$ for which $F$ can be represented by an $n$-mosaic.
\end{definition} 

By the construction, we have 

\begin{theorem}
The mosaic number $m(F)$ is an invariant for immersed surface-links. 
\end{theorem}

\begin{figure}[h!]
  \centering
  \includegraphics[width = 12cm]{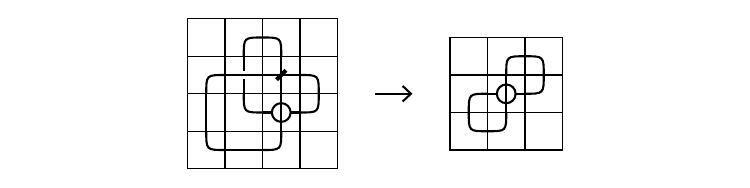}
  \caption{Two mosaic diagrams of $F_{1}$}\label{ExaMosaic1}
\end{figure}

\begin{figure}[h!]
  \centering
  \includegraphics[width = 12cm]{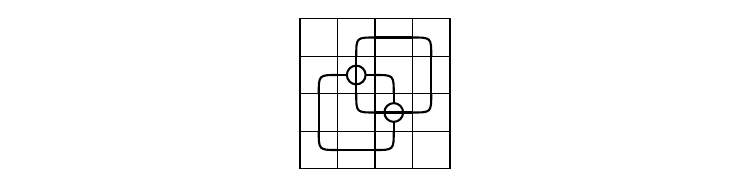}
  \caption{A mosaic diagram of $2_{1}^{0,0}(2)$}\label{ExaSimple}
\end{figure}

For example, the immersed surface-link $F_{1}$ whose diagram is described in Figure \ref{ExaSMGD} has a $4$-mosaic diagram and it can be changed to a $3$-mosaic diagram as shown in Figure \ref{ExaMosaic1}, by using mosaic moves.
Notice that when the mosaic number is $3$, there is exactly one case, and it is precisely $F_{1}$. 
It is also obvious that the mosaic number of $2_{1}^{0,0}(2)$ is $4$ as described in Figure \ref{ExaSimple}. 

The \textit{chs-index} of an immersed surface-link $F$ can be defined as the minimal number of $\mathrm{chs}(D)$ over all singular marked graph diagrams $D$ of $F$ where $\mathrm{chs}(D)$ is the total number of crossings, marked vertices, and singular vertices in $D$.

\begin{lemma}
Let $F$ be an immersed surface-link with the chs-index greater than or equal to $8$.
If $F$ has exactly one singular vertex, then $m(F)\geq 6$.
\end{lemma}

\begin{figure}[h!]
  \centering
  \includegraphics[width = 12cm]{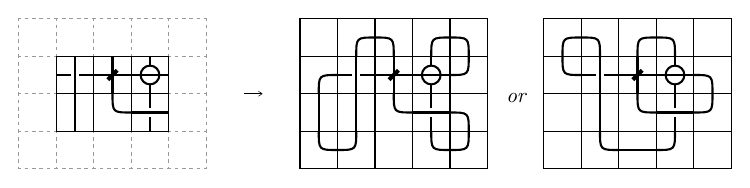}
  \caption{Twofold rule}\label{Twofold}
\end{figure}

The {\it inner tiles} are tiles of a $(m, n)$-mosaic diagram whose all connection points except boundary points meet the connection points of contiguous tiles. 
For example, see the $(2, 3)$-mosaic diagram in the left side of Figure \ref{Twofold}.
If we add boundary tiles to given inner tiles in order to connect the boundary connection points so that the whole tiles are suitably connected, then we get two mosaic diagrams. 
This is known as the {\it twofold rule} in \cite{OHLL} and it is useful to calculate the mosaic numbers.
Note that from the suitably connected condition, crossings, marked vertices, and singular vertices can be located in the position of inner tiles. 

\begin{proof}
When $\mathrm{chs}(F)\geq 10$, it is clear that $m(F)\geq 6$. 
Consider the cases that $\mathrm{chs}(F)=8$ or $9$. 
Suppose that $m(F)=5$. 
Let $D$ be a $5$-mosaic diagram of $F$. 
Since the number of inner tiles in $D$ is nine and since it has only one singular vertex, there exists at least one row such that it has no singular vertices, up to rotations, as depicted in Figure \ref{OneRow}.
\begin{figure}[h!]
  \centering
  \includegraphics[width = 12cm]{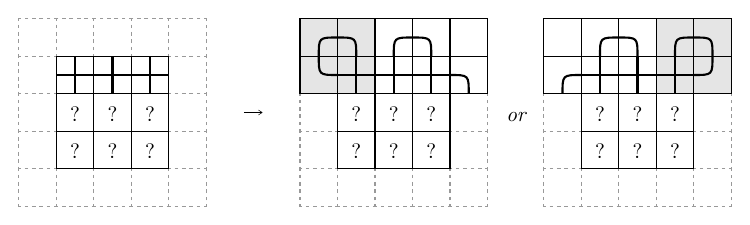}
  \caption{One row without singular vertices applying the twofold rule}\label{OneRow}
\end{figure}
By the twofold rule, the resulting diagrams have at least one kink with one classical crossing or marked vertex.
\begin{figure}[h!]
  \centering
  \includegraphics[width = 12cm]{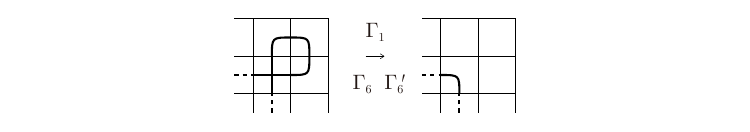}
  \caption{A kink with one crossing or marked vertex}\label{Kink}
\end{figure}
By applying Yoshikawa moves $\Gamma_{1}, \Gamma_{6},$ or $\Gamma_{6}'$, the kink can be removed in Figure \ref{Kink}. 
This contradicts that $\mathrm{chs}(F)=8$ or $9$.
Hence $m(F)\geq 6$.
\end{proof}

\section*{Acknowledgement}
All authors equally contribute this paper. 
The work of Seonmi Choi was supported by Basic Science Research Program
through the National Research Foundation of Korea (NRF) funded by the Ministry of Education (No. 2021R1I1A1A01049100) and the National Research Foundation of Korea (NRF) grant funded by the Korean government (MSIT) (No. 2022R1A5A1033624).
The work of Jieon Kim was supported by Young Researchers Program through the National Research Foundation of Korea (NRF), funded by the Ministry of Education, Science and Technology (NRF- 2018R1C1B6007021).

\bigskip

\noindent
\textsc{Nonlinear Dynamics and Mathematical Application Center \\
Kyungpook National University \\
Daegu, 41566, Republic of Korea}

\bigskip

\noindent
\textsc{Department of Mathematics\\ 
Pusan National University\\ 
Busan, 46241, Republic of Korea}

\end{document}